\documentclass{gtart_h}  

\def\ifplaintex{\expandafter\ifx\csname documentclass\endcsname\relax}

\ifplaintex 
\else
\fi

\expandafter\ifx\csname beginpicture\endcsname\relax
\expandafter\ifx\csname documentclass\endcsname\relax
\input pictex \else
\input prepictex \input pictex \input postpictex \fi\fi

\def\gt{{\mathsurround=0pt\it $\cal G\mskip-2mu$eometry \&\ 
$\cal T\!\!$opology}}        

\def\gtp{{\mathsurround=0pt\it $\cal G\mskip-2mu$eometry \&\ 
$\cal T\!\!$opology $\cal P\!$ublications}}  

\def\lognumber#1{\def\thelognumber{#1}}
\def\volumenumber#1{\def\thevolumenumber{#1}}
\def\papernumber#1{\def\thepapernumber{#1}}
\def\volumeyear#1{\def\thevolumeyear{#1}}

\def\pagenumbers#1#2{\def\startpage{#1}\def\finishpage{#2}}
\def\published#1{\def\publishdate{#1}}
\def\proposed#1{\def\theproposer{#1}}
\def\seconded#1{\def\theseconders{#1}}
\def\received#1{\def\receiveddate{#1}}

\def\accepted#1{\def\accepteddate{#1}}
\def\asciititle#1{\def\theasciititle{#1}}
\def\covertitle#1{\def\thecovertitle{#1}}
\def\coverauthors#1{\def\thecoverauthors{#1}}
\def\asciiauthors#1{\def\theasciiauthors{#1}}
\def\asciiaddress#1{\def\theasciiaddress{#1}}
\def\asciiemail#1{\def\theasciiemail{#1}}

\long\def\asciiabstract#1{\long\def\theasciiabstract{#1}}

\let\\\par\let\thelognumber\relax
\let\thevolumenumber\relax\let\thepapernumber\relax
\let\thevolumeyear\relax\let\thesamplenumber\relax\let\startpage\relax
\let\finishpage\relax\let\publishdate\relax\let\receiveddate\relax
\let\reviseddate\relax\let\accepteddate\relax\let\theasciititle\relax
\let\thecovertitle\relax\let\theasciiauthors\relax\let\theasciiaddress\relax
\let\theasciiabstract\relax
\let\theasciiemail\relax\let\theshortauthors\relax\let\theshorttitle\relax
\let\thecoverauthors\relax

\long\def\maketitlep{   

\count0=\startpage

\gt\hfill      
\beginpicture
\setcoordinatesystem units <0.33truein, 0.33truein> point at 2.2 0.9
\setplotsymbol ({$\cal G$})
\plotsymbolspacing=9truept
\circulararc 315 degrees from 0 1 center at 0 0
\setplotsymbol ({$\cal T$})
\circulararc 315 degrees from 1 -1 center at 1 0
\endpicture
\break
{\small\ifx\thesamplenumber\relax 
Volume \else Sample
\fi\thevolumenumber\ (\thevolumeyear)
\startpage--\finishpage\nl
Published: \publishdate}
\vglue 0.5truein plus 0.4fil minus 0.1truein

{\parskip=0pt\leftskip 0pt plus 1fil\def\\{\par\smallskip}{\ifplaintex\large
\else\Large\fi\bf\thetitle}\par\medskip}   

\vglue 0pt plus 0.1fil 

{\parskip=0pt\leftskip 0pt plus 1fil\def\\{\par}{\sc\theauthors}
\par\medskip}

\vglue 0pt plus 0.1fil 

{\small\parskip=0pt\let\newline\\
{\leftskip 0pt plus 1fil\def\\{\par}{\sl\theaddress}\par}
\expandafter\ifx\theemail\relax    
\relax\else\vglue 5pt plus 0.02fil minus 2pt\def\\{\stdspace{\rm 
and}\stdspace} 
\cl{Email:\stdspace\tt\theemail}\fi
\ifx\theurl\relax                  
\relax\else\vglue 5pt plus 0.02fil minus 2pt\def\\{\stdspace{\rm 
and}\stdspace}
\cl{URL:\stdspace\tt\theurl}\fi\par}

\vglue 7pt plus 0.3fil minus 3pt

{\bf Abstract}
\vglue 5pt plus 0.1fil minus 2pt

\theabstract

\vglue 7pt plus 0.3fil minus 3pt

{\bf AMS Classification numbers}\quad Primary:\quad \theprimaryclass

Secondary:\quad \thesecondaryclass

\vglue 5pt plus 0.3fil minus 2pt

{\bf Keywords:}\quad \thekeywords

\vglue 10pt plus 0.5fil minus 5pt

{\small  Proposed: \theproposer\hfill Received: \receiveddate\nl
Seconded: \theseconders\hfill 
\ifx\reviseddate\relax                         
Accepted: \accepteddate                        
\else
Revised: \reviseddate                          
\fi}
\eject
}       

\let\maketitlepage\maketitlep
\let\maketitle\maketitlepage

\font\phead=cmsl9 scaled 950
\font=cmsl9 scaled 1050
\font\pnum=cmbx10 scaled 913
\font\lnum=cmbx10 
\font\pfoot=cmsl9 scaled 950
\font=cmsl9 scaled 1050
\ifplaintex
\headline{\vbox to 0pt{\vskip -4.5mm\line{\small\phead\ifnum
\count0=\startpage ISSN 1364-0380 (on line)
1465-3060 (printed) \hfill {\pnum\folio}\else\ifodd\count0\def\\{ }%
\ifx\theshorttitle\relax\thetitle\else\theshorttitle\fi\hfill{\pnum\folio}
\else\def\\{ and }{\pnum\folio}\hfill\ifx\theshortauthors\relax\theauthors
\else\theshortauthors\fi\fi\fi}\vss}}
\footline{\vbox to 0pt{\vglue 0mm\line{\small\pfoot\ifnum\count0=\startpage
\copyright\ \gtp\hfill\else
\gt, Volume \thevolumenumber\ (\thevolumeyear)\hfill\fi}\vss
}}
\else
\makeatletter
\def\@oddhead{{\small\lhead\ifnum\count0=\startpage ISSN 1364-0380 (on line)
1465-3060 (printed) \hfill {\lnum\number\count0}\else\ifodd\count0
\def\\{ }\ifx\theshorttitle\relax \thetitle \else\theshorttitle\fi\hfill
{\lnum\number\count0}\else\def\\{ and }{\lnum\number\count0}
\hfill\ifx\theshortauthors\relax 
\theauthors\else\theshortauthors\fi\fi\fi}}\def\@evenhead{\@oddhead}
\def\@oddfoot{\small\lfoot\ifnum\count0=\startpage\copyright\ \gtp\hfill\else
\gt, Volume \thevolumenumber\ (\thevolumeyear)\hfill\fi}
\def\@evenfoot{\@oddfoot}
\makeatother
\fi

\newwrite\gtoutfile
\long\gdef\makeheadfile{  
{\def\\{, }\def\s{ }
\immediate\openout\gtoutfile head.xxx
\immediate\write\gtoutfile{Proxy-for: \ifx\theasciiauthors\relax
\theauthors\else\theasciiauthors\fi\s<\ifx\theasciiemail\relax\theemail\else\theasciiemail\fi>}
\immediate\write\gtoutfile{\noexpand\\}
\immediate\write\gtoutfile{Authors: \ifx\theasciiauthors\relax
\theauthors\else\theasciiauthors\fi}
{\def\\{ }\immediate\write\gtoutfile{Title: \ifx\theasciititle\relax
\thetitle\else\theasciititle\fi}}
\immediate\write\gtoutfile{Subj-class: GT or SG or MG etc}
\immediate\write\gtoutfile{MSC-class: \theprimaryclass\ifx\thesecondaryclass\relax\else, \thesecondaryclass\fi}
\immediate\write\gtoutfile{Journal-ref: Geom. Topol. \thevolumenumber
(\thevolumeyear) \startpage-\finishpage}
\immediate\write\gtoutfile{Comments: Published by Geometry and Topology at}
\immediate\write\gtoutfile{\s\s http://www.maths.warwick.ac.uk/gt/GTVol\thevolumenumber/paper\thepapernumber.abs.html}
\immediate\write\gtoutfile{\noexpand\\}
\immediate\write\gtoutfile{}
\ifx\theasciiabstract\relax
\immediate\write\gtoutfile{\theabstract}\else
\immediate\write\gtoutfile{\theasciiabstract}\fi
\immediate\write\gtoutfile{}
\immediate\write\gtoutfile{\noexpand\\}
\immediate\write\gtoutfile{}
\immediate\closeout\gtoutfile}}  

\def\maketitlepage{\maketitlep\makeheadfile}
\let\maketitle\maketitlepage

\lognumber{529}
\received{25 November 2004}
\volumenumber{9}\papernumber{9}\volumeyear{2005}
\pagenumbers{299}{313}   
\published{28 January 2005}
\accepted{24 January 2005}
\proposed{Yasha Eliashberg}
\seconded{Leonid Polterovich, David Gabai}

\usepackage{amssymb,amsmath}

\newtheorem{theoreme}{Th\'eor\`eme}[section]

\newtheorem{lemme}[theoreme]{Lemme}
\newtheorem{corollaire}[theoreme]{Corollaire}

\newcommand{\R}{\mathbb{R}}
\newcommand{\Z}{\mathbb{Z}}

\newcommand{\N}{\mathbb{N}}

\newcommand{\bdry}{\partial}
\newcommand{\s}{\vskip.1in}

\newcommand{\F}{\mathcal{F}}
\newcommand{\tr}{\textrm{tr}}

\newcommand{\be}{\begin{enumerate}}
\newcommand{\ee}{\end{enumerate}}

\begin{document}

\title{Homologie de contact des vari\'et\'es toro\"\i dales}
\asciititle{Homologie de contact des varietes toroidales}
\covertitle{Homologie de contact des vari\noexpand\'et\noexpand\'es toro\noexpand\"\noexpand\i dales}

\author{Fr\'ed\'eric Bourgeois\\Vincent Colin}
\coverauthors{Fr\noexpand\'ed\noexpand\'eric Bourgeois\\Vincent Colin}
\asciiauthors{Frederic Bourgeois\\Vincent Colin}

\address{Universit\'e Libre de Bruxelles, D\'epartement de 
Math\'ematiques CP 218\\Boulevard du Triomphe, 1050 Bruxelles, Belgium}

\secondaddress{Universit\'e de Nantes, Laboratoire de Math\'ematiques 
Jean Leray\\2 rue de la Houssini\`ere, BP 92208, 44322 Nantes Cedex 3, France}

\asciiaddress{Universite Libre de Bruxelles, Departement de
  Mathematiques CP 218 \\Boulevard du Triomphe, 1050 Bruxelles,
  Belgium\\and\\Universite 
de Nantes, Laboratoire de Mathematiques Jean Leray\\2 rue
de la Houssiniere, BP 92208, 44322 Nantes Cedex 3, France}

\asciiemail{fbourgeo@ulb.ac.be, Vincent.Colin@math.univ-nantes.fr}
\gtemail{\mailto{fbourgeo@ulb.ac.be}{\rm\qua 
and\qua}\mailto{Vincent.Colin@math.univ-nantes.fr}}

\begin{abstract}
We show that contact homology distinguishes infinitely 
many tight contact structures on any orientable, toroidal, 
irreducible 3--manifold. As a consequence of the contact 
homology computations, on a very large class of toroidal 
manifolds, all known examples of universally tight contact 
structures with nonvanishing torsion satisfy the Weinstein 
conjecture.
\smallskip

{\bf R\'esum\'e}

\smallskip
On montre que l'homologie de contact distingue
une infinit\'e de structures de contact tendues sur toute
vari\'et\'e toro\"\i dale irr\'eductible et orientable de dimension trois.
En cons\'equence des calculs d'homologie de contact,
sur une tr\`es large classe de vari\'et\'es toro\"\i dales,
tous les exemples de structures de contact universellement tendues
de torsion non nulle connus v\'erifient la conjecture de Weinstein.
\end{abstract}

\asciiabstract{
We show that contact homology distinguishes infinitely 
many tight contact structures on any orientable, toroidal, 
irreducible 3-manifold. As a consequence of the contact 
homology computations, on a very large class of toroidal 
manifolds, all known examples of universally tight contact 
structures with nonvanishing torsion satisfy the Weinstein 
conjecture.

Resume 

On montre que l'homologie de contact distingue
une infinite de structures de contact tendues sur toute
variete toroidale irreductible et orientable de dimension trois.
En consequence des calculs d'homologie de contact,
sur une tres large classe de varietes toroidales,
tous les exemples de structures de contact universellement tendues
de torsion non nulle connus verifient la conjecture de Weinstein.}



\maketitlepage

\section{Introduction}
Dans cet article, on prouve que l'homologie de contact distingue
une infinit\'e de structures de contact tendues sur toute
vari\'et\'e toro\"\i dale, irr\'eductible et orientable
de dimension trois.

\begin{theoreme}\label{theoreme : toroidal}
{\rm (Voir aussi \cite{Co1, Co2}, \cite{HKM})}\qua
Toute vari\'et\'e toro\"\i dale, irr\'eductible, orientable, close
de dimension trois porte
une infinit\'e de structures de contact hypertendues
deux \`a deux non isotopes et distingu\'ees par leur homologie de contact.
\end{theoreme}

Cette approche d\'emontre que l'homologie de contact \cite{EGH}
d\'etecte la torsion \cite{Gi, Co1} dans son cadre le plus g\'en\'eral.
La famille de structures de contact \'etudi\'ee
est la m\^eme que dans \cite{Co1,Co2} et \cite{HKM}.
Ici, on montre que tous ses
membres sont deux \`a deux non isotopes, alors
qu'avec les arguments de \cite{Co1,Co2,HKM}, on obtient
seulement l'existence d'une sous-famille infinie non explicite ayant
cette propri\'et\'e.
En revanche, on n'obtient pas en g\'en\'eral que les structures
consid\'er\'ees sont deux \`a deux non isomorphes.

Dans le cas des fibr\'es en tores sur le cercle, l'homologie 
de contact est toujours de rang infini. On d\'efinit alors un 
invariant num\'erique en enrichissant l'homologie de contact au 
moyen de la fonctionnelle d'action pour distinguer une infinit\'e de 
structures de contact sur ces vari\'et\'es.

Les exemples \'etudi\'es sont les seules structures
de contact universellement tendues connues~-- et conjecturalement les
seules~-- sur les vari\'et\'es toro\"\i dales irr\'e\-duc\-tibles. On
montre ici (th\'eor\`eme~\ref{theoreme : Weinstein}) que,
parmi celles-ci, sur
``presque toutes les vari\'et\'es", toutes celles de torsion non
nulle ont une homologie de contact non nulle. En particulier,
elles v\'erifient la conjecture de Weinstein: leurs champs de
Reeb poss\`edent tous une orbite p\'eriodique. Par ailleurs, dans
ce cas, les structures non universellement tendues sont
virtuellement vrill\'ees, et v\'erifient \'egalement, d'apr\`es
Hofer~\cite{Ho}, la conjecture de Weinstein. D'apr\`es \cite{CGH},
les structures de contact tendues de torsion nulle 
sur une vari\'et\'e donn\'ee sont en nombre
fini; d'apr\`es ce qui pr\'ec\`ede, ce sont potentiellement les
seules pour lesquelles la conjecture de Weinstein reste ouverte.
Une preuve de la conjecture de Weinstein pour les structures
de contact port\'ees par un livre ouvert planaire
a r\'ecemment \'et\'e annonc\'ee par Abbas, Cieliebak et
Hofer \cite{ACH}. Il est vraisemblable que les structures
de contact tendues de torsion non nulle soient toutes
de genre sup\'erieur \`a $1$, auquel cas nos exemples
seraient disjoints de ceux de \cite{ACH}.

Une structure de contact $\xi$ est {\it hypertendue} si elle
poss\`ede un champ de Reeb $R$ sans orbite p\'eriodique
contractible. Dans le cas o\`u la vari\'et\'e $V$ est toro\"\i
dale, toutes les orbites de $R$ sont d'ordre infini dans le groupe
fondamental $\pi_1 (V)$. En particulier, tout rappel de $\xi$ dans
un rev\^etement fini de $V$ est hypertendu, et donc tendu
d'apr\`es Hofer \cite{Ho}. Lorsque $V$ est toro\"\i dale, $\pi_1
(V)$ est r\'esiduellement fini, et la structure $\xi$
universellement tendue.

Toutes les vari\'et\'es consid\'er\'ees par la suite sont orient\'ees,
et les structures de contact positives pour cette orientation.
Les bords sont orient\'es par la r\`egle de la ``normale
sortante en premier". Si $A$ est une vari\'et\'e orient\'ee,
$-A$ d\'esigne la vari\'et\'e $A$ munie de l'orientation oppos\'ee.

{\bf Remerciements}\qua VC remercie le CNRS pour son accueil en d\'el\'egation
au cours de l'ann\'ee 2003--2004.

\section{Construction d'une famille infinie}

Dans cette section, on d\'ecrit la construction des structures de contact
non isotopes dans le th\'eor\`eme \ref{theoreme : toroidal}.
Pour plus de simplicit\'e, on se restreint au cas
o\`u la torsion appara\^\i t dans une seule classe d'isotopie
de tore. L'\'enonc\'e g\'en\'eral est rejet\'e ˆ la section \ref{section : general}.

On part d'une vari\'et\'e $V$ irr\'eductible, toro\"\i dale, close et
orient\'ee de dimension trois.
Soit $T\subset V$ un tore incompressible et $T^2 \times [a,b]$ un voisinage
tubulaire de $T\simeq T^2 \times \{ (a+b)/2\}$.
On note $V'$ la vari\'et\'e $V'=V\setminus T^2 \times ]a,b[$.
On se donne une surface essentielle minimale coorient\'ee $S$ dans $V'$
(i.e.\ incompressible, bord-incompressible, sans disque parall\`ele au bord, 
sans composante close, et qui minimise
la norme de Thurston dans $H_2 (V',\bdry V' , \Z )$),
qui rencontre les deux composantes $T^2 \times \{ a, b\}$ de
$\partial V'$.

On a besoin de mesurer des pentes dans $T^2 \times [a,b+2k\pi ]$, $k\in \Z$.
Pour cela, on choisit une base $([l],[m])$
de $H_1 (T^2,\Z )$, et on consid\`ere la projection
$$
p \co  T^2 \times [a,b+2k\pi ] \rightarrow T^2 = (\R /\Z )^2 :  ((x,y),t) \rightarrow (x,y)
$$ 
de fibre $[a,b+2k\pi ]$.
\`A toute courbe non contractible $\gamma$ dans $T^2 \times [a,b+2k\pi ]$,
on associe sa {\it pente} \'egale \`a
$(p,q)/\sqrt{p^2 +q^2} \in S^1 \subset \R^2$, o\`u $p$ et $q$
sont d\'efinis par l'identit\'e
$$
[p(\gamma )]= p[l]+q[m] \in H_1 (V,\Z ) .
$$
Comme $\gamma$ est orient\'ee, elle est aussi coorient\'ee.
Son {\it c\^one positif} $C(\gamma )$ est le demi-cercle ferm\'e
($\subset S^1$), d\'elimit\'e par les pentes
de $\gamma$ et $-\gamma$, et contenant les pentes
des courbes $\delta$ de $T^2$ qui ont une intersection
homologique $p(\gamma) \cdot \delta$ positive.

Une fois ces notations en place, on se repose sur la construction 
suivante \cite{Co3}:

\begin{theoreme}\label{theoreme : tore} Pour tout $\epsilon >0$,
il existe une forme de contact
hypertendue $\alpha_\epsilon$ sur $V$ et deux r\'eels $c,c'>0$
avec les propri\'et\'es suivantes:
\begin{itemize}
\item le champ de Reeb $R_\epsilon$ de $\alpha_\epsilon$ est positivement 
transversal \`a $S$, tangent \`a $T^2 \times \{ a,b\}$;
\item $R_\epsilon$ vaut $\cos t\partial_x -\sin t\partial_y$
sur $T^2 \times [a+c ,b-c']$ ($\alpha_\epsilon$ vaut $\cos tdx-\sin tdy$);
\item sur $T^2 \times [a,a+c ]$ et $T^2 \times [b-c' ,b]$,
$R_\epsilon$ est positivement transversal \`a un feuilletage en anneau
$\F_a \times [a,a+c ]\cup (-\F_b )\times [b-c' ,b]$, o\`u $\F_a$
et $\F_b$ sont des feuilletages en cercles de $T^2 \times \{ a\}$ et $T^2 \times \{ b\}$
dont les pentes sont $\epsilon$--proches de celles d'une composante
de, respectivement, $\partial S \cap T^2 \times \{ a\}$ et $\partial S
\cap T^2 \times \{ b\}$.
\end{itemize}
\end{theoreme}

On peut en particulier trouver une
forme de contact hypertendue $\alpha_{\epsilon }$ qui satisfait
aux conclusions du th\'eor\`eme~\ref{theoreme : tore} avec de surcro\^\i t
$(b-a)/(2\pi )\leq 1$.
Pour tout $k\in \N$, on note alors $\alpha_{\epsilon ,k}$ la forme de contact
hypertendue obtenue en rempla\c cant dans $V$
le produit de contact $(T^2 \times [a+c, b-c' ],\cos tdx- \sin tdy)$
par $(T^2 \times [a+c, b-c' +2k\pi ],\cos tdx- \sin tdy)$.
Les formes $\alpha_{\epsilon ,k}$ sont hypertendues, donc les
structures $\xi_{\epsilon ,k}=\ker \alpha_{\epsilon ,k}$ 
sont universellement tendues.

On montre que, pour $\epsilon$ assez petit fix\'e,
les structures $\xi_{\epsilon ,k}$ sont, dans la plupart des cas,
deux \`a deux non isotopes.

\section{Preuve du th\'eor\`eme $1$.$1$}

La preuve utilise l'homologie de contact, d\'efinie dans la section \ref{section : hc},
et est subdivis\'ee de la mani\`ere suivante.  Dans la section \ref{section : cas1}, 
on \'etudie le cas o\`u la pente de $\partial S \cap T^2 \times \{ a\}$ est diff\'erente 
de celle de $\partial S\cap T^2 \times \{ b\}$, dans la section \ref{section : cas2} 
le cas o\`u ces pentes sont identiques, et dans la section \ref{section : cas3}
le cas des fibr\'es en tores sur le cercle et des vari\'et\'es obtenues par collage 
de deux fibr\'es tordus en intervalles sur la bouteille de Klein.

\subsection{Homologie de contact}  \label{section : hc}

L'homologie de contact est un invariant des structures de contact introduit par
Eliashberg, Givental et Hofer \cite{EGH}. Nous utiliserons ici la version la plus 
simple de cet invariant : l'homologie de contact cylindrique.

Pour la d\'efinir, on fixe une classe d'homotopie libre
de lacets $h$ de $V$, et on consid\`ere le $\Z$--module $C^h_*$ librement 
engendr\'e par les orbites p\'eriodiques de $R_{\epsilon ,k}$ dans la classe $h$. 
Si la forme de contact $\alpha_{\epsilon, k}$ est choisie de mani\`ere
g\'en\'erique, ces orbites p\'eriodiques sont non d\'eg\'en\'er\'ees, et
on peut associer \`a chaque orbite $\gamma$ son indice de 
Conley--Zehnder $\mu(\gamma)$. La graduation de $\gamma$ dans
le module $C^h_*$ est donn\'ee par $\mu(\gamma) - 1$.

On munit la structure de contact $\xi_{\epsilon,k}$ d'une structure complexe
$J$ compatible avec la forme symplectique $d\alpha_{\epsilon ,k}$. On
\'etend ensuite $J$ en une structure presque complexe, toujours not\'ee $J$,
sur  la symplectisation $(\R \times V, d(e^t \alpha_{\epsilon ,k}))$ de 
$(V, \xi_{\epsilon ,k})$, en posant $J \frac{\partial}{\partial t} = R_{\epsilon ,k}$.
On consid\`ere les espaces de modules des cylindres $J$--holomorphes 
dans $(\R \times V, d(e^t \alpha_{\epsilon ,k}))$,
convergeant lorsque $t \to \pm\infty$ vers des orbites p\'eriodiques de 
$R_{\epsilon ,k}$ dans la classe $h$. 

Lorsque la classe $h$ ne contient que des orbites simples de 
$R_{\epsilon ,k}$, on peut choisir $J$ g\'en\'erique\-ment pour que 
ces espaces de modules soient des vari\'et\'es lisses \cite{D}.
Comme les formes $\alpha_{\epsilon ,k}$ sont hypertendues, 
leurs champs de Reeb $R_{\epsilon ,k}$ ne comptent aucune
orbite p\'eriodique contractible, et donc les cylindres $J$--holomorphes
ne peuvent d\'eg\'en\'erer en faisant appara\^itre un plan complexe.
Par cons\'equent, les espaces de modules consid\'er\'es sont compacts 
\cite{BEHWZ}. En particulier, l'homologie de contact cylindrique est
bien d\'efinie dans cette situation.

On d\'efinit une diff\'erentielle lin\'eaire $d \co  C^h_* \rightarrow C^h_{*-1}$ 
au moyen des cylindres $J$--holomorphes 
dans la symplectisation $(\R \times V, d(e^t \alpha_{\epsilon ,k}))$, 
de sorte que
le coefficient de $\gamma_2$ dans $d\gamma_1$ compte le nombre de cylindres 
$J$--holomorphes rigides convergeant vers $\gamma_1$ lorsque $t \to +\infty$
et vers $\gamma_2$ lorsque $t \to -\infty$. Les espaces de modules sont
orient\'es  \cite{BM} de mani\`ere coh\'erente, de mani\`ere \`a pouvoir compter
les cylindres $J$--holomorphes rigides avec un signe.

Cette diff\'erentielle $d$ satisfait $d^2 = 0$, en vertu d'un th\'eor\`eme de collage
\cite{Bo2} tr\`es semblable \`a celui utilis\'e pour l'homologie de Floer. On a 
ainsi d\'efini un complexe $(C^h_*, d)$. L'homologie $HC^h_* (V,\xi ,\Z )$ de ce 
complexe, appel\'ee homologie de contact cylindrique, est un invariant 
de la classe de conjugaison de $\xi$ par tout diff\'eomorphisme qui fixe $h$, 
et en particulier par isotopie.

On montre que pour une classe $h$ bien choisie, les seules
orbites p\'eriodiques dans $h$ proviennent du produit
$T^2 \times [a+c,b-c' +2k\pi]$. On en d\'eduit un
calcul explicite de $HC_*^h (V,\xi_{\epsilon ,k} ,\Z )$
(corollaire~\ref{corollaire : homologie}).

Dans les sections \ref{section : cas1} et \ref{section : cas2}, 
la vari\'et\'e $V$ n'est pas un fibr\'e en tores sur le cercle
ni obtenue par collage de deux fibr\'es tordus en intervalles 
sur la bouteille de Klein. L'\'etude de ces deux cas est rejet\'ee
dans la section \ref{section : cas3}.

\subsection{La pente de $\partial S \cap T^2 \times \{ a\}$
est diff\'erente de celle de $\partial S\cap T^2 \times \{ b\}$}  \label{section : cas1}

Dans ce cas, on peut choisir $\epsilon$ assez petit de sorte que
le secteur
$$H=S^1 \setminus (C(\partial S \cap T^2 \times \{ a\} )
\cup C(-(\partial S \cap T^2 \times \{ b\} ))\cup C(\F_a )
\cup C(-\F_b ))$$ soit de cardinal infini.

Dans les lemmes qui suivent, on oriente les composantes
de bord de tout anneau $A$ dans le m\^eme sens de sorte qu'elles soient
librement homotopes dans $A$.

\begin{lemme}\label{lemme : anneau} On suppose qu'aucune composante de
$V'$ n'est diff\'eomorphe
\`a $T^2 \times [0,1]$ ni \`a un fibr\'e tordu en intervalles
sur la bouteille de Klein. Il existe des
pentes $p_a$  et $p_b$ telles que, pour tout anneau
$A$ immerg\'e dans $V'$, $\partial A\subset \partial V'$, les pentes des
deux composantes de $\partial A \subset T^2 \times [a,b]$ appartiennent
toutes deux \`a $\{\pm p_a ,\pm p_b \}$, ou soient identiques.
\end{lemme}

\begin{proof}[D\'emonstration] Si $A$ n'est pas essentiel, il est homotope
\`a un anneau inclus dans $\partial V'$ et les pentes
de ses deux composantes de bord sont identiques.

Sinon, $A$ est homotope \`a un anneau $A'$ inclus dans une
composante $V''$ de la d\'ecomposition de Jaco--Shalen--Johanson de
$V'$ qui est un fibr\'e de Seifert (d'apr\`es le {\it Mapping
Theorem} de Jaco--Shalen \cite{JS}). Les fibres de $V''$
rencontrent $T^2 \times \{ a\}$ ou $T^2 \times \{ b\}$; on note
$p_a$ ou $p_b$ leur pente dans $T^2 \times [a,b]$.

Par hypoth\`eses, $V''$ n'est pas un fibr\'e tordu en
intervalles sur une bouteille
de Klein ni le tore \'epais ni un tore solide.
D'apr\`es Jaco--Shalen (Lemme II.2.8, \cite{JS}),
les pentes des composantes de bord de $A$ sont
$\pm p_a$ o\`u $\pm p_b$.
\end{proof}

Soit $\gamma_0$ une courbe
de $T^2 \times [a,b]$ dont la pente est dans $H \setminus \{ \pm p_a
,\pm p_b \}$.
On note $h$ sa classe d'homotopie libre.

\begin{lemme}\label{lemme : homotopie} Si aucune composante de $V'$
n'est un tore \'epaissi ou un fibr\'e tordu en intervalles
sur la bouteille de Klein,
toute orbite de $R_{\epsilon ,k}$ dans la classe $h$ est incluse dans
$$
T^2 \times [a+c,b-c' +2k\pi ]
$$ 
et est de m\^eme pente que $\gamma_0$.
\end{lemme}

\begin{proof}[D\'emonstration]
Soit $\gamma_1$ une orbite p\'eriodique de $R_{\epsilon ,k}$
librement homotope \`a $\gamma_0$. Elle est incluse soit
dans $T^2 \times [a+c,b-c' +2k\pi ]$, soit dans $T^2 \times ([a,a+c]
\cup [b-c'+2k\pi, b+2k\pi ])$, soit dans $V'$.

L'existence d'une homotopie entre $\gamma_0$ et $\gamma_1$
fournit l'existence d'un anneau immerg\'e $A=\phi (S^1 \times [0,1])$,
$\partial A=\gamma_0 \cup \gamma_1$.

Par une homotopie de $A$, on se ram\`ene \`a un anneau encore lisse, toujours
not\'e $A$, mais qui est \'egalement transversal \`a $T^2 \times \{ a,b+2k\pi \}$
et tel que toutes ses intersections avec $T^2 \times \{ a,b+2k\pi \}$
soient non contractibles dans $V$. Pour obtenir cette deuxi\`eme
propri\'et\'e, on conjugue le fait que les tores $T^2 \times \{ a, b+2k\pi \}$
sont incompressibles avec le r\'esultat d'un algorithme classique de Haken.
A la source dans $S^1 \times [0,1]$, les pr\'eimages par $\phi$ des courbes
d'intersection forment une collection de courbes disjointes et parall\`eles.
Elles d\'ecoupent $S^1 \times [0,1]$ en une collection d'anneaux.
On note $\delta_1$,...,$\delta_{p-1}$ leur image par $\phi$ et
$A_1$,...,$A_p$ celle des anneaux.
Par convention, $\gamma_0 ,\delta_1 \subset A_1$,
$\gamma_1 ,\delta_{p-1} \subset A_p$ et, pour $2\leq k\leq p-1$,
$\partial A_k =\delta_{k-1} \cup \delta_{k}$.

Par construction, $A_1 \subset T^2 \times [a,b+2k\pi ]$ et donc la pente
de $\delta_1$ est identique \`a celle de $\gamma_0$.
D'apr\`es le lemme~\ref{lemme : anneau},
la pente de $\delta_2$ est soit la m\^eme que celle de $\delta_1$,
soit $\delta_1$ et $\delta_2$ ont pour pente $\pm p_a$ et $\pm p_b$.
Par hypoth\`ese, ce n'est pas le cas pour $\gamma_0$, donc
ce n'est \'egalement pas le cas pour $\delta_1$.
On en d\'eduit que la pente de $\delta_2$ est la m\^eme que celle
de $\delta_1$. Une r\'ecurrence imm\'ediate montre que, pour $1\leq p-1$,
le pente de $\delta_k$ est la m\^eme que celle de $\gamma_0$.

Si $\gamma_1 \subset T^2 \times ([a,a+c] \cup [b-c'+2k\pi ,b+2k\pi
])$, alors on a aussi 
$$
A_p \subset T^2 \times ([a,a+c] \cup
[b-c'+2k\pi ,b+2k\pi ]) .
$$ 
La pente de $\gamma_1$ est dans ce cas
\'egale \`a celle de $\gamma_0$. Ceci met en contradiction le
choix de $h$ et le fait que le champ de Reeb soit positivement
transversal \`a un feuilletage $\F_a \times [a,a+c] \cup (-\F_b )
\times [b-c'+2k\pi ,b+2k\pi ]$.

Dans le cas o\`u $\gamma_1 \subset V'$,
on observe que $[\delta_{p-1} ]\in H_1 (V',\Z )$ a, dans $V'$, une intersection
n\'egative avec $[S] \in H_2 (V',\partial V',\Z )$.
Comme $[\gamma_1 ]=[\delta_{p-1} ]$, l'intersection de $[\gamma_1 ]$
avec $[S]$ est \'egalement n\'egative. Ceci contredit le fait que
toutes les intersections de $\gamma_1$ avec $S$ sont positives
(car le champ de Reeb est positivement transversal \`a $S$).

Pour finir, si $\gamma_1$ est dans $T^2 \times [a+c,b-c'+2k\pi ]$,
sa pente est la m\^eme que celle de $\gamma_0$.
\end{proof}

\begin{corollaire}\label{corollaire : homologie}
Pour $\epsilon$ assez petit, si $V$ n'est pas un fibr\'e en tores sur le
cercle et si aucune composante de $V'$ n'est un fibr\'e tordu en intervalles
sur la bouteille de Klein, $HC_0^h (V,\xi_{\epsilon ,k} )\simeq \Z^k$.
\end{corollaire}

\begin{proof}[D\'emonstration] On a exactement $k$ cercles d'orbites 
de $R_{\epsilon ,k}$
dans la classe $h$ qui sont toutes situ\'ees dans $T^2 \times
[a+c,b-c'+2k\pi ]$.
Toutes les orbites de ces familles ont la m\^eme p\'eriode.
Par cons\'equent, les seuls cylindres holomorphes entre des orbites 
de la classe $h$
dans la symplectisation de $(V,\alpha_{\epsilon ,k} 
)$ sont les cylindres verticaux.
En s'inspirant de techniques \`a la 
Morse--Bott \cite{Bo1},
on peut perturber la forme de contact $\alpha_{\epsilon,k}$ au moyen 
d'une fonction de
Morse $f$ ayant exactement un maximum et un minimum 
sur chaque cercle d'orbites
de la classe $h$, de sorte que chaque 
cercle donne lieu \`a deux orbites non d\'eg\'en\'er\'ees, 
l'une de degr\'e $-1$ et l'autre de degr\'e $0$.
Ces orbites sont les 
g\'en\'erateurs du complexe de cha\^\i ne pour l'homologie de 
contact.
Apr\`es la perturbation, les seuls cylindres holomorphes 
d'index $1$ correspondent aux
trajectoires du gradient de $f$. Leurs 
contributions \`a la diff\'erentielle du complexe de
cha\^\i ne 
s'annulent donc deux \`a deux, et
chaque g\'en\'erateur du complexe 
de cha\^\i ne
donne lieu \`a un g\'en\'erateur de l'homologie de contact.
\end{proof}

Dans le cas o\`u une seule composante $V_1$ de $V'$ (attach\'ee
par exemple le long de $T^2 \times \{ b\}$) est un fibr\'e tordu
en intervalles sur la bouteille de Klein, on se ram\`ene
au cas pr\'ec\'edent par passage au rev\^etement. 
Soit en effet $\pi \co  T^2 \times [-1,1 ]\rightarrow V_1$ une application
de rev\^etement de degr\'e $2$
qui est compatible avec la forme de contact $\cos tdx -\sin tdy$
de $T^2 \times [-1,1]$. Si $T^2 =(\R /\Z )^2$, la projection
$\pi$ est le passage au quotient pour la relation d'\'equivalence:
$(x,y,t) \sim (x,-y,-t)$.

On munit $V_1$ de la structure de contact quotient.
On d\'esigne par $V_2$ l'autre composante de $V'$.

On note $(W,\eta_{\epsilon ,k} )$ le rev\^etement de degr\'e $2$ de
$(V,\xi_{\epsilon ,k} )$ obtenu en collant
deux copies du compl\'ementaire de l'int\'erieur de $V_1$ dans $V$
au rev\^etement $\pi^{-1} (V_1 )$. Dans $W$, le tore \'epais 
$$
P=T^2 \times [a,b+2k\pi ]\cup \pi^{-1} (V_1 )
\cup T^2 \times [a,b+2k\pi ]
$$ 
joue le r\^ole de
$T^2 \times [a,b+2k\pi ]$ dans $V$: la structure $\eta_{\epsilon ,k}$
y est conjugu\'ee \`a un produit $(T^2 \times [a',b'] ,\ker
(\cos tdx-\sin tdy))$, $4k\pi \leq b'-a'\leq (4k+2)\pi$.
Le compl\'ementaire de $P$ dans $W$ est constitu\'e de deux copies
de $V_2$
munies de la structure de contact $\ker \alpha_{\epsilon ,0}$.

D'apr\`es le th\'eor\`eme~\ref{theoreme : tore}, pour tout
$\epsilon >0$, il existe $c,c'>0$ et une forme de contact
$\beta_{\epsilon , k}$ pour $\eta_{\epsilon ,k}$ qui vaut
$\alpha_{\epsilon ,0}$ sur $W\setminus P$, $\cos tdx-\sin tdy$
sur $T^2 \times [a+c,b+2k\pi ]\cup \pi^{-1} (V_1 )
\cup T^2 \times [a,b-c'+2k\pi ]$
et dont le champ de Reeb est positivement transversal \`a un feuilletage
$\F_a \times [a,a+c]\cup (-\F_b  )\times [b-c',b]$ dont les deux
pentes dans $P$ sont $\epsilon$--proches de
celles des copies de $\partial S$ dans $\partial P$.

De plus, les pentes dans $P$ des deux copies
du bord de $S$ situ\'ees de part et d'autre
de $P$ sont oppos\'ees.

Gr\^ace \`a l'\'etude pr\'ec\'edente,
on montre que, pour $\epsilon$ assez petit, les structures $\eta_{\epsilon ,k}$
sont deux \`a deux non isotopes dans $W$.
C'est donc aussi le cas des structures de la famille $(\xi_{\epsilon,k} )_{k\in \N}$.

\subsection{La pente de $\partial S \cap T^2 \times \{ a\}$
est \'egale \`a celle de $\partial S\cap T^2 \times \{ b\}$}  \label{section : cas2}

Si $T$ s\'epare $V$, la surface $S$ se d\'ecompose elle
aussi en deux composantes situ\'ees de part et d'autre
de $T^2\times [a,b]$.
Il suffit de changer l'orientation d'une des deux composantes
de $S$ pour se ramener au cas pr\'ec\'edent: on d\'efinit
une nouvelle famille de structures
$(\xi'_{\epsilon ,k} )_{k\in \N}$ deux \`a deux non isotopes.

Si $T$ ne s\'epare pas $V$, on modifie la structure $\xi_{\epsilon ,k}$
sur $V$ en collant \`a la restriction de $\alpha_{\epsilon ,0}$ sur $V\setminus
T^2 \times ]a+c,b-c'[$
le produit  $(T^2 \times [a+c, b-c'+(2k+1)\pi ], \cos tdx-\sin tdy)$.
La structure $\xi'_{\epsilon ,k}$ ainsi obtenue n'est pas orientable.
Elle le devient si on consid\`ere son rappel dans
le rev\^etement $W$ de degr\'e $2$ de $V$ obtenu en
collant deux copies de $V\setminus T^2 \times ]a+c,b-c'[$ \`a deux copies de
$T^2 \times [a+c, b-c'+(2k+1)\pi ]$.
Dans $W$, on r\'ecup\`ere une forme de contact $\beta_{\epsilon ,k}$
en prenant $\cos tdx-\sin tdy$ sur $\pi^{-1}
(T^2 \times [a+c,b-c'+(2k+1)\pi ])$ et, respectivement, $\alpha_{\epsilon ,0}$
et $-\alpha_{\epsilon ,0}$ sur les deux exemplaires de
$V\setminus T^2 \times ]a+c,b-c'+(2k+1)\pi [$.

On se fixe une composante $P$ de $\pi^{-1} (T^2 \times [a+c,
b-c'+(2k+1)\pi ])$ dans $W$. L'autre est not\'ee $Q$. Le champ de
Reeb associ\'e \`a la forme $\beta_{\epsilon ,k}$ est transversal
aux deux rappels $S_1$ et $S_2$ de $S$, positivement dans un cas
($S_1$) et n\'egativement dans l'autre. On consid\`ere alors la
surface $S'=S_1 \cup S_2$. On choisit $\epsilon$ petit pour que
$H=S^1 \setminus (C(\partial S' )\cup C(\F_a ) \cup C(-\F_b ))$
soit de cardinal infini. Comme pr\'ec\'edemment, on note $p_a$ et
$p_b$ les pentes des \'eventuelles composantes fibr\'ees
adjacentes \`a $P$ dans la d\'ecomposition de
Jaco--Shalen--Johanson du compl\'ementaire de $P$ dans $W$. 
On fixe une classe homotopie libre $h$ d'une courbe de $P$ de pente 
dans $H\setminus \{ \pm p_a ,\pm p_b \}$.

On montre, comme dans la partie pr\'ec\'edente, que toute orbite
du flot de Reeb de $\beta_{\epsilon ,k}$ dans la classe $h$ doit provenir de
l'une des pr\'eimages de $T^2 \times [a+c,b-c'+(2k+1)\pi ]$
(dans le cas qui nous int\'eresse, aucune composante
de $V'$ n'est un tore \'epaissi, un fibr\'e tordu en
intervalles sur une bouteille de Klein, ou un tore solide).
De plus, les seules orbites de $Q$ homotopes \`a une orbite de
$P$ sont homotopes \`a des orbites de pente $\pm p_a$ ou $\pm p_b$.
En cons\'equence, on obtient le corollaire:

\begin{corollaire} Pour $\epsilon >0$ assez petit, $HC_0^h (W,\ker
\beta_{\epsilon ,k} ,\Z )\simeq \Z^k$.
\end{corollaire}

On en d\'eduit imm\'ediatement que, pour
$\epsilon$ fix\'e assez petit, les structures
de la famille $(\xi'_{\epsilon ,k} )_{k\in \N}$ sont deux \`a deux non isotopes
sur $V$.

\subsection{Deux cas particuliers}   \label{section : cas3}

On traite ici le cas des fibr\'es en tores sur le cercle
et des vari\'et\'es obtenues par collage de deux
fibr\'es tordus en intervalles sur la bouteille de Klein.
Dans ces deux cas, les familles de structures sont explicites.

Comme une structure de contact sur un fibr\'e tordu en intervalles sur la
bouteille de Klein est le quotient par $\Z_2$ d'une structure de
contact sur un fibr\'e en tores sur le cercle, ces
deux situations
sont analogues. Nous ne traiterons donc que le cas o\`u
$V = (\R/\Z)^2 \times [0,1]/\sim$, avec les identifications $(x,y,1) \sim
(A(x,y),0)$ pour $A \in SL(2,\Z)$. Les structures de contact
$\xi_{\epsilon,k} = \ker \alpha_{\epsilon,k}$ introduites
dans la section 2 sont isotopes aux structures de contact $\xi_k = \ker
\alpha_k$ d\'efinies par
$$
\alpha_k = \cos ((a+2k\pi)z+b) dx + \sin((a+2k\pi)z+b) dy .
$$
Soit $A^t$ la transpos\'ee de la matrice $A$.
Les nombres $0 < a, b \le 2\pi$ sont
choisis de sorte que $(\cos b,\sin b) =
A^t(\cos(a+b),\sin(a+b))$.

On distingue trois cas selon la monodromie $A$:
\begin{enumerate}
\item[(i)] Si $| \tr A | < 2$, la monodromie
$A$ est d'ordre fini dans $SL(2,\Z)$. Par cons\'equent,
toute classe
d'homotopie libre $h$ dans $T^2$ ne contient qu'un nombre fini
d'orbites de Reeb ferm\'ees.
\item[(ii)] Si $| \tr A | = 2$, la
monodromie $A$ est conjugu\'ee \`a
$$
\pm \left( \begin{array}{cc}
1 & n
\\
0 & 1
\end{array} \right)
$$
dans $SL(2,\Z)$ pour un certain
entier $n$. En particulier, il existe une classe d'homotopie libre
$h$ dans $T^2$ ne contenant qu'un nombre fini et non nul d'orbites
de Reeb ferm\'ees.
\item[(iii)] Si $| \tr A | > 2$, la monodromie ne
pr\'eserve aucune direction de pente rationnelle. Par
cons\'equent,
toute classe d'homotopie libre $h$ dans $T^2$ contient  un nombre
infini d'orbites de Reeb ferm\'ees, lorsque $k \ge 1$.
\end{enumerate}

Pour calculer la diff\'erentielle du complexe pour l'homologie de 
contact, nous aurons
besoin du lemme suivant.

\begin{lemme}
Soit $F \co  (\R \times S^1,j) \to (\R \times V,J)$
un cylindre holomorphe dans la symplectisation de $(V, \alpha_k)$
et aboutissant en $t = -\infty$ \`a une orbite de Reeb ferm\'ee dans 
$T^2 \times \{ z_0 \}$. 
Alors le cylindre $F(\R \times S^1)$ est contenu dans $\R 
\times T^2 \times \{ z_0 \}$.
\end{lemme}

\begin{proof}[D\'emonstration]
Supposons par 
contradiction qu'il existe $p \in \R \times S^1$ tel que 
$F(p) \notin \R \times T^2 
\times \{ z_0 \}$. Sans perte de g\'en\'eralit\'e, 
on peut supposer 
que $z \circ F(p) = z_0-\epsilon$ ou $z_0 + \epsilon$,
pour un petit $\epsilon > 0$. Pour fixer les id\'ees, on choisit
$z \circ F(p) = z_0 + \epsilon$.

Par le th\'eor\`eme de Sard, on peut choisir deux valeurs
r\'eguli\`eres $z_1$ et $z_2$ de la fonction $z \circ F$ sur
$\R \times S^1$, satisfaisant $z_0  < z_1 < z_2 < z_0 +\epsilon$.
Comme il n'y a pas d'orbite de Reeb ferm\'ee de classe $h$ dans
$T^2 \times [z_1, z_2]$ pour $\epsilon > 0$ suffisamment petit,
le sous-ensemble 
$C_{z_1,z_2} =  F^{-1}(\R \times T^2 \times [z_1,z_2])$ du
cylindre $\R \times S^1$ est  non vide, compact  et  bord\'e 
de lacets simples et lisses. 
De plus, les images de ces lacets dans $\R \times V$ sont soit 
contractibles, soit de classe $h$.

Par le th\'eor\`eme de Stokes, on a
$$
\int_{C_{z_1,z_2}} F^* d\alpha_k =
\int_{F^{-1}(\R \times T^2 \times \{ z_2 \}) } F^* \alpha_k 
- \int_{F^{-1}(\R \times T^2 \times \{ z_1 \}) } F^* \alpha_k  .
$$
Ces deux int\'egrales sont ais\'ement calculables, puisque la $1$--forme
$\alpha_k$ est ferm\'ee sur les sous-vari\'et\'es $T^2 \times \{ z_1 
\}$ et
$T^2 \times \{ z_2 \}$.
Leur diff\'erence est nulle ou du signe de $\cos((a+2k\pi)(z_2 - z_0))
- \cos((a+2k\pi)(z_1 - z_0)) < 0$ selon que la classe d'homologie de
$F^{-1}(\R \times T^2 \times \{ z_1 \})$ est nulle ou 
correspond \`a $h$.

D'autre part, comme $F$ est holomorphe, on doit avoir 
$F^* d\alpha_k \ge 0$ sur $\R \times S^1$. Par cons\'equent, 
la $2$--forme $d\alpha_k$ doit s'annuler 
sur $C_{z_1,z_2}$, de sorte que $F$ est un cylindre vertical 
au-dessus d'une orbite de Reeb. On obtient 
ainsi une contradiction.
\end{proof}

Comme dans le corollaire \ref{corollaire : homologie}, 
on d\'eduit de ce lemme que la 
diff\'erentielle du complexe 
d\'efinissant l'homologie de contact
$HC^h_0(V,\xi_k)$ est 
identiquement nulle.
Par cons\'equent, si $h$ est une classe 
d'homotopie libre dans $V$ contenant un nombre 
fini et non nul 
d'orbites de Reeb ferm\'ees, alors le rang de $HC^h_0(V,\xi_k)$ 
d\'epend de $k$. Plus pr\'ecis\'ement, si $m$ d\'esigne le plus
petit entier strictement positif tel que $A^m h = h$, alors
$$
\textrm{rang } HC^h_0(V,\xi_{k+1}) = \textrm{rang } HC^h_0(V,\xi_k) + m .
$$
En particulier, dans les cas (i) et (ii), les 
structures de contact $\xi_k$ sont non isotopes 
pour des valeurs 
diff\'erentes de $k$.

Si $h$ est une classe d'homotopie libre dans 
$V$ contenant un nombre infini 
d'orbites de Reeb ferm\'ees, 
l'homologie de contact est de rang infini pour tout $k$
et ne permet 
pas de conclure. On utilise alors l'homologie de contact en 
conjonction avec la
fonctionnelle d'action, qui associe \`a toute 
orbite de Reeb ferm\'ee sa p\'eriode. Soit
$\alpha$ une forme de 
contact pour une structure de contact $\xi$ sur $V$; on 
note
$HC^{h,T}_*(V,\alpha)$ l'ensemble des classes d'homologie dans 
$HC^h_*(V,\alpha)$ 
ayant un repr\'esentant de p\'eriode inf\'erieure 
ou \'egale \`a $T$. Notons que ceci d\'epend du
choix de la forme de 
contact $\alpha$ pour $\xi$. On d\'efinit
$$
I^h_*(V,\alpha) = 
\lim_{T \to \infty} \frac{\textrm{rang } HC^{h,T}_*(V,\alpha)}{\ln T} 
$$
lorsque la limite existe. Dans le cas contraire, il est encore 
possible de d\'efinir 
$\underline{I}^h_*(V,\alpha)$ et 
$\overline{I}^h_*(V,\alpha)$ en utilisant une limite 
inf\'erieure ou 
sup\'erieure respectivement.

\begin{lemme}
La quantit\'e 
$I^h_*(V,\alpha)$ ne d\'epend que de la classe d'isotopie de 
$\xi$.
\end{lemme}

\begin{proof}[D\'emonstration]
Consid\'erons la 
forme de contact $c\alpha$ pour $c > 0$. On a 
\begin{eqnarray*}
I^h_*(V,c\alpha) &=& \lim_{T \to \infty} 
\frac{\textrm{rang } HC^{h,T}_*(V,c\alpha)}{\ln T} \\
&=& \lim_{T \to 
\infty} \frac{\textrm{rang } HC^{h,\frac{T}{c}}_*(V,\alpha)}{\ln T} 
\\
&=& \lim_{T \to \infty} \frac{\textrm{rang } 
HC^{h,\frac{T}{c}}_*(V,\alpha)}{\ln \frac{T}{c}} 
\frac{\ln T - \ln 
c}{\ln T} \\
&=&  I^h_*(V,\alpha).
\end{eqnarray*}
Soient $\alpha_1$ 
et $\alpha_2$ deux formes de contact sur $V$. On dira que
$\alpha_1 < 
\alpha_2$ si on peut munir $\R \times V$ d'une forme symplectique 
$\omega$
exacte telle que $\omega =  d(e^t \alpha_2)$ lorsque $t \le 
0$ et $\omega =  d(e^{t-t_0} \alpha_1)$ 
lorsque $t \ge t_0$, pour un 
certain $t_0 > 0$. Si $\alpha_1$ et $\alpha_2$
sont des formes de 
contact pour des structures de contact isotopes sur $V$ alors il 
existe $c > 0$
tel que $c\alpha_1 < \alpha_2$. 

Si $\alpha_1 < 
\alpha_2$ alors $\textrm{rang } HC^{h,T}_*(V,\alpha_2) \le 
\textrm{rang } HC^{h,T}_*(V,\alpha_1)$. En effet, en vertu des 
propri\'et\'es fonctorielles de l'homologie
de contact \cite{EGH}, au 
cobordisme symplectique $(\R \times V, \omega)$
correspond un 
isomorphisme 
$$
\Phi \co  HC^h_*(V,\alpha_2) \to HC^h_*(V,\alpha_1)
$$ 
tel que
$\Phi(HC^{h,T}_*(V,\alpha_2)) \subset 
HC^{h,T}_*(V,\alpha_1)$.
Par cons\'equent, si $\alpha_1 < \alpha_2$ 
alors $I^h_*(V,\alpha_2) \le I^h_*(V,\alpha_1)$. 

La conclusion suit 
de ces deux propri\'et\'es de $I^h_*$.
\end{proof} 

En vertu de ce r\'esultat, nous utiliserons plut\^ot la notation 
$I^h_*(V,\xi)$, qui ne fait plus appara\^itre la forme de contact $\alpha$. 

Calculons maintenant $I^h_0(V,\xi_k)$ lorsque $|\tr A| > 2$, en 
consid\'erant plut\^ot la
vari\'et\'e de contact $(T^2 \times [0,1], 
\xi_k)$.
La classe d'homotopie libre $h$ dans $V$ correspond aux 
classes d'homotopie libres
$A^n h$ dans $T^2 \times [0,1]$, pour $n 
\in \Z$. Il y a $k$ ou $k+1$ orbites de Reeb
ferm\'ees dans la classe 
d'homotopie libre $A^n h$; cette alternative ne d\'epend que de 
$\xi_0$ lorsque $|n|$ est suffisamment grand. L'action de ces 
orbites est donn\'ee par 
$\| A^n h \| \simeq C \lambda^{|n|}$ 
lorsque $|n|$ est grand, o\`u $\lambda > 1$ est l'une 
des valeurs 
propres
de $A$. On en d\'eduit que $I^h_0(V,\xi_k)$ est \'egal \`a 
$\frac{2k}{\ln \lambda}$
ou $\frac{2(k+1)}{\ln \lambda}$. En 
particulier, les structures de contact $\xi_k$ ne sont pas
isotopes 
pour diff\'erentes valeurs de $k \ge 0$.

\section{G\'en\'eralisations}  \label{section : general}

L'\'etude propos\'ee ci-dessus s'adapte dans le cas o\`u la torsion
est non nulle dans plusieurs classes d'isotopies.

Plus pr\'ecisemment, soient $T_0$,...,$T_n$ des tores incompressibles
deux \`a deux disjoints et non parall\`eles dans
une vari\'et\'e close et irr\'eductible $V$.
On note $T_i \times ]a_i ,b_i [$ un voisinage tubulaire
de $T_i$ dans $V$ (choisis deux \`a deux disjoints)
et $S$ une surface minimale essentielle coorient\'ee
dans $V'=V\setminus \cup_i (T_i \times ]a_i ,b_i [)$
qui rencontre toutes les composantes de $\partial V'$.

Comme dans le th\'eor\`eme~\ref{theoreme : tore}, on munit
$V$, pour $\epsilon >0$, d'une forme de contact $\alpha_\epsilon$,
dont le champ de Reeb $R_\epsilon$
poss\`ede, pour tout $0\leq i\leq n$, les propri\'et\'es suivantes:
\begin{itemize}
\item $R_\epsilon$ est positivement transversal \`a $S$, tangent
\`a $T_i \times \{ a_i ,b_i \}$;
\item $R_\epsilon$ vaut $\cos t\partial_x -\sin t\partial_y$
sur $T_i \times [a_i +c ,b_i-c']$ ($\alpha_\epsilon$ vaut $\cos tdx-\sin tdy$);
\item sur $T_i \times [a_i ,a_i +c ]$ et $T_i \times [b_i -c' ,b_i ]$,
$R_\epsilon$ est positivement transversal \`a un feuilletage en anneau
$\F_{a_i} \times [a_i ,a_i +c ]\cup (-\F_{b_i} )\times [b_i -c' ,b_i ]$,
o\`u $\F_{a_i}$ et $\F_{b_i}$ sont des
feuilletages en cercles de $T_i \times \{ a_i \}$ et $T_i \times \{ b_i \}$
dont les pentes sont $\epsilon$--proches de celles d'une composante
de, respectivement, $\partial S \cap T_i \times \{ a_i \}$ et $\partial S
\cap T_i \times \{ b_i \}$.
\end{itemize}

Sur chaque produit $T_i \times ]a_i ,b_i [$,
on remplace le produit de contact $(T_i \times [a_i +c, b_i -c' ],\cos tdx- \sin tdy)$
par $(T_i \times [a_i +c, b_i -c' +2k_i \pi ],\cos tdx- \sin tdy)$,
pour un certain $k_i \in \N$.
On obtient ainsi une certaine forme de contact
$\alpha_{\epsilon ,k_0 ,...,k_n}$.

En reprenant le m\^eme raisonnement que pr\'ec\'edemment, et avec
les notations introduites ci-dessus, on montre le th\'eor\`eme suivant:

\begin{theoreme}\label{theoreme : Weinstein} Si, pour un certain $0\leq i\leq n$,
$S \cap T_i \times \{ a_i \}$ et $S \cap T_i \times \{ b_i \}$
ne sont pas de pente oppos\'ee (en projection dans $T_i$) et si $k_i \geq 1$,
alors il existe une classe d'homotopie $h_i$ de courbe
dans $T_i$ pour laquelle:
$$HC_*^{h_i} (V, \xi_{k_0 ,...,k_n} ,\Z )\simeq \Z^{k_i}.$$
En particulier, la structure $\xi_{k_0 ,...,k_n}$ v\'erifie
la conjecture de Weinstein.
\end{theoreme}

Tous les exemples de structures de contact universellement
tendues connus sont de ce type. La restriction principale
dans l'\'enonc\'e du th\'eor\`eme~\ref{theoreme : Weinstein} est
topologique et porte sur la pente de
$\partial S \cap T_i \times \{ a_i, b_i \}$.

\end{document}